\documentclass[a4paper,12pt,reqno]{amsart}
\usepackage[english,francais]{babel}
\usepackage{amsfonts}
\usepackage{amsmath}
\usepackage{amssymb}
\usepackage{mathrsfs}
\usepackage{nomencl}
\usepackage{delarray}
\usepackage{amsthm}
\usepackage{blkarray}
\DeclareMathOperator{\Sp}{Sp}

\usepackage[colorlinks]{hyperref}

\numberwithin{equation}{section}
\theoremstyle{plain}

\newtheorem{thm}{Theorem}[section]

\newtheorem{prop}[thm]{Proposition}

\theoremstyle{definition}
\newtheorem{defn}[thm]{Definition}

\newtheorem{rem}[thm]{Remark}
\newtheorem{ex}[thm]{Example}

\def\RR{{\mathbb R}}

\def\C{{\mathbb C}}
\def\Gh{{\widehat{G}}}

\def\H{\mathcal{H}}

\def\VN{{\rm VN}}
\def\SU2{{\rm SU(2)}}

\begin{document}
\title[Fourier multipliers and group von Neumann algebras]
{Fourier multipliers and group von Neumann algebras}

\author[Rauan Akylzhanov]{Rauan Akylzhanov}
\address{Rauan Akylzhanov:
  \endgraf
  Department of Mathematics
  \endgraf
  Imperial College London
  \endgraf
  180 Queen's Gate, London SW7 2AZ
  \endgraf
  United Kingdom
  \endgraf
  {\it E-mail address} {\rm r.akylzhanov14@imperial.ac.uk}
}

\author[Michael Ruzhansky]{Michael Ruzhansky}
\address{
  Michael Ruzhansky:
  \endgraf
  Department of Mathematics
  \endgraf
  Imperial College London
  \endgraf
  180 Queen's Gate, London SW7 2AZ
  \endgraf
  United Kingdom
  \endgraf
  {\it E-mail address} {\rm m.ruzhansky@imperial.ac.uk}
  }

\selectlanguage{english}%
\thanks{The second author was supported by the
Leverhulme Research Grant RPG-2014-02 and by
 the EPSRC Grant EP/K039407/1.  
 No new data was collected or generated during the course of the research.}
\date{\today}
  \medskip


\maketitle
\begin{abstract}
\selectlanguage{english}%
In this paper we establish the $L^p$-$L^q$ boundedness of Fourier multipliers
	on locally compact separable unimodular groups for
	the range of indices $1<p\leq 2 \leq q<\infty$. Our approach is based on the operator algebras techniques.
	The result depends on a 
	version of the Hausdorff-Young-Paley inequality that we establish on 
	general locally compact separable unimodular groups.
	In particular, the obtained result implies the corresponding 
	H\"ormander's Fourier multiplier theorem 
	on $\mathbb{R}^{n}$ and the corresponding known results for Fourier 
	multipliers on compact Lie groups.

{\it To cite this article: 
R. Akylzhanov, M. Ruzhansky, 
C. R. Acad. Sci. Paris, Ser. *** (20**).}\vskip 0.5\baselineskip
	
\selectlanguage{francais}%
	\noindent{\bf R\'esum\'e} \vskip 0.5\baselineskip \noindent
\textbf{Multiplicateurs de Fourier et alg\'ebres de von Neumann.}
Dans cet article nous \'etablissons des $L^p$-$L^q$ bornitudes de multiplicateurs de Fourier sur les groupes unimodulaire localement compacts pour $1<p\leq 2 \leq q<\infty$.
Notre approche est bas\'ee sur la technique des alg\'ebres des op\'erateurs.
Pour cela nous prouvons une version de l'in\'egalit\'e de Hausdorff-Young sur les groupes unimodulaires localement compacts. En particulier, le r\'esultat obtenu implique le th\'eor\`eme de H\"ormander sur les multiplicateurs de Fourier 
dans $ \mathbb{R}^{n} $ et des r\'esultats d\'ej\`a connus associ\'es aux multiplicateurs de Fourier sur les groupes de Lie compacts.

{\it Mot-cl\'es: alg\'ebres d'op\'erateurs, groupes unimodulaires localement compacts, op\'erateurs pseudo-diff\'erentiels.}

{\it Pour citer cet article: R. Akylzhanov, M. Ruzhansky, 
C. R. Acad. Sci. Paris, Ser. *** (20**).}
\end{abstract}


\selectlanguage{english}

\selectlanguage{francais}
\section*{Version fran\c{c}aise abr\'eg\'ee}
L'objectif de cet article est de donner des conditions suffisantes pour des $L^p$-$L^q$ bornitudes de multiplicateurs de Fourier sur les groupes unimodulaires localement compacts.

 La $L^p$-$L^q$ bornitude des multiplicateurs de Fourier \`a \'et\'e r\'ecemment \'etudi\'ee dans le contexte des groupes compacts (\cite{ANR2015}). Le r\'esultat de ce papier (Th\'eor\`eme \ref{THM:upper-bound}) g\'en\'eralise des r\'esultats connus: 
H\"ormander \cite[Th\'eor\'eme 1.11]{Hormander:invariant-LP-Acta-1960} et 
 \cite[Th\'eor\'eme 3.1]{ANR2015}, maintenant
pour les groupes localement compacts unimodulaires s\'eparables $G$.
L'hypoth\`ese de s\'eparablilit\'e et d'unimodularit\'e pour les groupes localement compacts 
 peut \'etre consid\'er\'ee comme naturelle car permettant d'utiliser des r\'esultats de base de l'analyse de Fourier de von Neumann comme, par example, la formule de Plancherel  (Segal \cite{Segal1950}). 
Pour une discussion plus d\'etaill\'ee des op\'erateurs pseudo-diff\'erentiels dans ce context 
nous nous r\'ef\'erons \`a \cite{MR2015}, mais nous notons que dans cet article nous n'avons pas besoin que la groupe sois, par example, de type I ou II. Certains r\'esultats sur les multiplicateurs $L^p$ de Fourier dans l'esprit du th\'eor\'eme de Mikhlin-H\"ormander sont \'egalement connus sur les groupes: par example Coifman et Weiss \cite{coifman+weiss_lnm} pour le groupe SU(2), \cite{RuWi2015} pour les groupes compacts generaux et \cite{Fischer-Ruzhansky:FM-graded} pour les  groupes de Lie gradu\'es. Leur approche des multiplicateurs de Fourier $L^p$ est diff\'erente de celle propos\'ee dans cet article, ce qui nous permet d'enlever l'hypoth\'ese que la groupe sois compact ou nilpotent.

\selectlanguage{english}

\section{Introduction}
The aim of this paper is to give sufficient conditions for the $L^p$-$L^q$ boundedness of Fourier multipliers on  locally compact separable unimodular groups $G$.
To put this 
in context, we recall that in \cite[Theorem 1.11]{Hormander:invariant-LP-Acta-1960}, 
Lars H\" ormander established sufficient conditions on the symbol for the $L^p$-$L^q$ boundedness of Fourier multipliers on $\RR^n$:
\begin{equation}
\label{EQ:Hormander-estimate}
\|A\|_{L^p(\RR^n)\to L^q(\RR^n)}
\lesssim
\sup_{\substack{s>0}}s
\left(\,\,
\int\limits_{\substack{\xi\in\RR^n\\ |\sigma_A(\xi)|\geq s}}d\xi
\right)^{\frac1p-\frac1q},\quad 1<p\leq 2 \leq q <+\infty.
\end{equation}
Here, as usual, the Fourier multiplier $A$ acts by multiplication on the Fourier transform side,
i.e. $\widehat{Af}(\xi)=\sigma_{A}(\xi)\widehat{f}(\xi)$, $\xi\in\RR^{n}.$

The $L^{p}$-$L^{q}$ boundedness of Fourier multipliers has been also recently
investigated in the context of compact Lie groups in \cite{ANR2015}.
The result of this paper (Theorem \ref{THM:upper-bound})
generalises known results: H\"ormander's inequality \eqref{EQ:Hormander-estimate} as well as \cite[Theorem 3.1]{ANR2015} 
to the setting of general locally compact separable unimodular groups $G$.
See also \cite{ANR2016} for a wide collection of results in the particular case of the group SU(2).

The assumption for the locally compact group to be separable and unimodular may be viewed as
quite natural allowing one to use basic results of von Neumann-type 
Fourier analysis, such as, for example,
Plancherel formula (see Segal \cite{Segal1950}).
For a more detailed discussion of pseudo-differential operators in such settings
we refer to \cite{MR2015}, but we note that in this paper we do not need to assume that
the group is, for example, of type I or type II.
Some results on $L^{p}$-Fourier multipliers in the spirit of H\"ormander-Mikhlin theorem
are also known on groups. See, for example, Coifman and Weiss
\cite{coifman+weiss_lnm} for the case of the group SU(2),
\cite{RuWi2015} for general compact Lie groups, and
\cite{Fischer-Ruzhansky:FM-graded} for graded Lie groups.
The approach to the $L^{p}$-Fourier multipliers 
is different from the technique proposed in this paper
allowing us to avoid here making an assumption that the group is compact or nilpotent.
Still, the nilpotent setting as in the book \cite{Fischer-Ruzhansky:book} allows for the use of symbols, something that we do not rely upon in the general setting in this note.

The proofs of the results of this note and possible extensions to non-unimodular groups or to non-invariant operators can be found in \cite{AR2016}.

\section{Notation and preliminaries}
Let $M\subset \mathcal{L}(\H)$ be a semifinite von Neumann algebra acting in a Hilbert space $\H$ with a trace $\tau$.
By an operator (in $S(M)$) we shall always mean a $\tau$-measurable operator affiliated with $M$ in the sense of Segal \cite{Segal1953}. We also refer to Dixmier \cite{VNA-Dixmier-1981} for the general background for the constructions below.

\begin{defn}
Let $A$ be an operator with the polar decomposition $A=U|A|$ and let
$|A|=\int\limits_{\Sp(|A|)}\lambda dE_{\lambda}$ be the spectral decomposition of $|A|$.
For each $t>0$, we define the generalised $t$-th singular numbers by
\begin{equation}
\label{EQ:mu-t}
\mu_t(A):=\inf\{\lambda\geq 0 \colon \tau(E_{\lambda}(|A|))\leq t\}.
\end{equation}
\end{defn}
These  $\mu_t(A)$ were first introduced by Fack as $t$-th generalised singular values of $A$, see \cite{FackKosaki1986} (and also below) for definition and properties.

As a noncommutative extension (\cite{Kosaki1981}) of the classical Lorentz spaces, we define Lorentz spaces $L^{p,q}(M)$ associated with a semifinite von Neumann algebra as follows:
\begin{defn} 
\label{DEF:Lorenz-spaces}
For $1\leq p <\infty$, $1\leq q \leq \infty$, denote by $L^{p,q}(M)$ the set of all operators $A\in S(M)$ satisfying
\begin{equation}
\|A\|_{L^{p,q}(M)}
:=
\left(
\int\limits^{+\infty}_0 
\left(
t^{\frac1p}\mu_t(A)
\right)^q
\frac{dt}{t}
\right)^{\frac1q}<+\infty.
\end{equation}
We define $L^{p}(M):=L^{p,p}(M)$, and
for $q=\infty$, we define
\begin{equation}
\|A\|_{L^{p,\infty}(M)}
:=
\sup_{t>0}t^{\frac1p}\mu_t(A).
\end{equation}
\end{defn}
\medskip
Let now $G$ be a locally compact unimodular separable group and let $M$  be the group von Neumann algebra ${\VN}_{L}(G)$ generated by the left regular representations $\pi_L(g)\colon f(\cdot)\to f(g^{-1}\cdot)$ with $g\in G$.

For $f\in L^1(G)\cap L^2(G)$, we say that $f$ on $G$ {\em has a Fourier transform} 
whenever the convolution operator 
\begin{equation}
\label{EQ:L_f}
L_fh(x):=(f\ast h)(x)=\int\limits_{G}f(g)h(g^{-1}x)\,dg
\end{equation}
 is a $\tau$-measurable operator with respect to $\VN_L(G)$.

\section{Results}

Our first result is a version of the Hausdorff-Young-Paley inequality on
locally compact unimodular separable groups. 

\begin{thm}[Hausdorff-Young-Paley inequality]
	\label{THM:HYP-LCG}
	Let $G$ be a locally compact unimodular separable group. 
	Let $1<p\leq b \leq p'<\infty$. 
	If a positive function $\varphi(t)$, $t\in\RR_+$, satisfies condition 
	 \begin{equation}
 \label{EQ:weak_symbol_estimate2}
	 M_\varphi:=\sup_{s>0}s\int\limits_{\substack{t\in\RR_+\\ \varphi(t)\geq s }}dt<\infty,
	 \end{equation}
then we have
	\begin{equation}
	\label{EQ:HYP-LCG}
	\left(
	\int\limits_{\RR_+}
	\left(
	\mu_t(L_f)
	{\varphi(t)}^{\frac1b-\frac{1}{p'}}
	\right)^b
	dt
	\right)^{\frac1b}
	\lesssim
	M_\varphi^{\frac1b-\frac1{p'}}
	\|f\|_{L^p(G)}.
	\end{equation}
\end{thm}
This inequality, besides being of interest on its own, is crucial in proving the multiplier
Theorem \ref{THM:upper-bound}. For $b=p'$, it gives Kunze's Hausdorff-Young inequality, see
e.g. \cite{Kunze1958}, and for $b=p$ it can be viewed as an analogue of the Paley inequality.

\begin{defn}
\label{DEF:FM}
 A linear operator $A$  is said to be a right {Fourier multiplier on $G$} if it is affiliated and $\tau$-measurable with respect to $\VN_L(G)$.
\end{defn}
It can be readily seen that $A$ being affiliated with $\VN_L(G)$ means that $A$ commutes with right translations on $G$, i.e. that $A$ is right-invariant. Thus, the right Fourier multipliers on $G$ (in the sense to Definition \ref{DEF:FM}) are precisely the right-invariant operators on $G$ that are
$\tau$-measurable with respect to $\VN_L(G)$.

In the following statements, to unite the formulations,
we adopt the convention that the sum or the integral over an
empty set is zero, and that $0^{0}=0$.

\begin{thm}
\label{THM:upper-bound}
 Let $1<p\leq 2 \leq q<+\infty$ and suppose that $A$ is a right Fourier multiplier on a locally compact separable unimodular  group $G$. Then we have
\begin{equation}
\label{EQ:LCG-FM-upper}
\|A\|_{L^p(G)\to L^q(G)}
\lesssim
\sup_{s>0}
s
\left[
\int\limits_{\substack{t\in\RR_{+}\colon \mu_t(A)\geq s}}dt
\right]^{\frac1p-\frac1q}.
\end{equation}
For $p=q=2$ inequality \eqref{EQ:LCG-FM-upper} is sharp, i.e.
\begin{equation}
\label{EQ:LCG-FM-upper:sharpness}
\|A\|_{L^2(G)\to L^2(G)}
=
\sup_{t\in\RR_+}\mu_t(A).
\end{equation}
Using the noncommutative Lorentz spaces $L^{r,\infty}$ 
with  $\frac1r=\frac1p-\frac1q$, $p\not=q$,
we can also write 
\eqref{EQ:LCG-FM-upper} as
\begin{equation}
\|A\|_{L^p(G)\to L^q(G)}
\lesssim
\|A\|_{L^{r,\infty}(\VN_L(G))}.
\end{equation}
\end{thm}

\begin{rem}
\label{REM:mu-t}
Let $G$ be a compact Lie group and let $\mathcal L$ be a self-adjoint positive left invariant operator on $G$ with discrete spectrum $\{\lambda_{k}\}\subset\mathbb R$. Let $\varphi:[0,\infty)\to (0,\infty)$ be a strictly decreasing strictly positive function. Then it can be shown that Theorem \ref{THM:upper-bound} implies that $\varphi(L)$ is bounded from $L^{p}(G)$ to $L^{q}(G)$, $1<p\leq 2 \leq q<+\infty$, provided that
\begin{equation}\label{EQ:phiL}
\sup_{k>0} k^{\frac1p-\frac1q} \varphi(\lambda_{k})<\infty. 
\end{equation}
For example, let us fix a sub-Riemannian structure  $(H,g)$ on $G$, i.e.
we choose a subbundle $H$ of the tangent bundle $TG$ of $G$ satisfying the H\"ormander bracket condition and a smooth metric $g$  on $H$.
The Carnot-Caratheodory metric induced by $g$ gives rise to the so-called Hausdorff measure on $G$. The intrinsic sub-Laplacian $-\Delta_b$ associated with the Hausdorff measure on $G$ can be shown to be a self-adjoint operator and has discrete spectrum (see e.g. \cite{HK2015}).
Let us denote by $\lambda_1(g),\lambda_2(g),\ldots,\lambda_k(g),\ldots$ the eigenvalues of the square root $\sqrt{1-\Delta_b}$ of $1-\Delta_b$.
The eigenvalue counting function $N_g(\lambda)$ of $\sqrt{1-\Delta_b}$ is given by
\begin{equation}
N_g(\lambda)=\sum\limits_{\substack{k\colon \lambda_k(g)\leq	 \lambda}}1
\end{equation}
(see \cite{Metivier1976} or \cite{HK2015}) and has the following asymptotics
\begin{equation}
N_g(\lambda)\cong \lambda^{Q}\quad \text{ as $\lambda\to+\infty$},
\end{equation}
where $Q$ is the Hausdorff dimension of $G$.
As a consequence, we obtain the asymptotic behaviour of the eigenvalues
\begin{equation}
\lambda_k(g)\cong k^{\frac1{Q}}, \quad \text{ as $k\to\infty$}.
\end{equation}
Then by Theorem \ref{THM:upper-bound} and taking $\varphi(\xi)=(1+\xi)^{-\frac{s}2}, \mathcal{L}=-\Delta_b$, the operator 
 $\varphi(-\mathcal{L})=(I-\Delta_b)^{-\frac{s}2}$ is bounded from $L^p(G)$ to $L^q(G)$, $1<p\leq 2 \leq q<+\infty$, provided that  
\begin{equation}
s\geq Q\left(\frac1p-\frac1q\right).
\end{equation}
If $H=TG$ and $g$ is the Killing form on $TG$, then $Q=n$ and $\Delta_b$ is the Laplacian on $G$ and we recover  the usual Sobolev embedding theorem in this case. 
\end{rem}

More details, extensions, and the proof of Remark \ref{REM:mu-t}  can be found in \cite{AR2016}.

In the case of locally compact abelian groups the statement of
Theorem \ref{THM:upper-bound} reduces nicely to (commutative) Lorentz spaces
defined in terms of non-increasing rearrangements of functions:

\begin{ex} Let $G$ be locally compact abelian group. The dual $\Gh$ of $G$ consists of continuous homomorphisms $\chi\colon G\to \C$. Then the group von Neumann algebra $\VN_L(G)$ is isometrically isomorphic to the multiplication algebra $L^{\infty}(\Gh)$, and
the operator $A$ acting by $\widehat{Af}(\chi)=\sigma_{A}(\chi)\widehat{f}(\chi)$ satisfies
 \begin{equation}
 \|A\|_{L^p(G)\to L^q(G)}
 \lesssim
 \sup_{s>0}
 s
 \left(
 \int\limits_{\substack{\chi\in\Gh\\ |\sigma_A(\chi)|\geq s}}
 d\chi
 \right)^{\frac1p-\frac1q}
 =
 \|\sigma_A\|_{L^{r,\infty}(\Gh)},
 \end{equation}
 where 
 \begin{equation}
 \|\sigma_A\|_{L^{r,\infty}(\Gh)}
 =
 \sup_{t>0}
 t^{\frac1r}
\sigma_A^*(t),\quad \frac1r=\frac1p-\frac1q, \; p\not=q.
 \end{equation}
Here $\sigma^*_A(t)$ is a non-increasing rearrangement of  the symbol $\sigma_A\colon\Gh\to\C$.
\end{ex}

In particular, for $G=\RR^n$, we recover H\"ormander's theorem 
(\cite[p. 106, Theorem 1.11]{Hormander:invariant-LP-Acta-1960}) for our range of $p$ and $q$:

\begin{rem} 
\label{REM:AR-H}
As a special case of $G=\RR^n$, Theorem \ref{THM:upper-bound} implies the H\"ormander multiplier estimate, namely,
\begin{equation}
\label{EQ:Hormander-estimate-2}
\|A\|_{L^p(\RR^n)\to L^q(\RR^n)}
\lesssim
\|\sigma_{A}\|_{L^{r,\infty}(\RR^n)}
=
\sup_{\substack{s>0}}s
\left(\,\,
\int\limits_{\substack{\xi\in\RR^n\\ |\sigma_A(\xi)|\geq s}}d\xi
\right)^{\frac1p-\frac1q},
\end{equation}
for the range $1<p\leq 2 \leq q <+\infty$, that was 
established in \cite[p. 106, Theorem 1.11]{Hormander:invariant-LP-Acta-1960}.
In this case it can be shown that 
\begin{equation}
\|A\|_{L^{r,\infty}(\VN_L(\RR^n))}
=
\|\sigma_A\|_{L^{r,\infty}(\RR^n)}.
\end{equation}
\end{rem}

We note that the same statements remain true also for left Fourier multipliers with the following modifications:

\begin{rem} 
\label{REM:AR-right}
The statement of Theorem \ref{THM:upper-bound} is also true for left Fourier multipliers,
i.e. for left-invariant operators on $G$ that are $\tau$-measurable with respect to 
the right group von Neumann algebra ${\VN}_{R}(G)$ generated by the right 
regular representations $\pi_R(g)\colon f(\cdot)\to f(\cdot \, g)$, $g\in G$.
In this case we have to replace all instances of ${\VN}_{L}(G)$ by ${\VN}_{R}(G)$.
\end{rem}

The $L^{p}$-$L^{q}$ boundedness of Fourier multipliers in the context of compact Lie groups
can be expressed in terms of their global matrix symbols: the $\tau$-measurability assumption is automatically satisfied under conditions of the following theorem, and so left Fourier multipliers 
$A$ are simply left-invariant operators on $G$; they act by
$\widehat{Af}(\xi)=\sigma_{A}(\xi)\widehat{f}(\xi)$ for $\xi\in\widehat{G}$,
$\sigma_{A}(\xi)\in\mathcal L(\mathcal H_{\xi})\simeq {\mathbb C}^{d_{\xi}\times d_{\xi}}$, where $d_{\xi}$ is the dimension of
the representation space $\mathcal H_{\xi}$ of $\xi\in\widehat G$. We refer to \cite{ANR2015} for details, and to \cite{RT}
or \cite{Ruzhansky+Turunen-IMRN} for the general theory of pseudo-differential operators 
on compact Lie groups.

\begin{thm}\label{THM:cmp}
Let $1<p\leq 2 \leq q <\infty$ and suppose that $A$ is a (left) 
Fourier multiplier on the compact Lie group $G$. Then we have
\begin{equation}
\label{EQ:CG-FM-upper}
\|A\|_{L^p(G)\to L^q(G)}
\lesssim
\sup_{\substack{s\geq 0}}
s
\left(
\sum\limits_{\substack{\xi\in\Gh\colon \|\sigma_A(\xi)\|_{\mathcal L(\mathcal H_{\xi})}\geq s}}
d^{2}_{\xi}
\right)^{\frac1p-\frac1q}.
\end{equation}
\end{thm}

One advantage of the estimate \eqref{EQ:CG-FM-upper}
(as well as of \eqref{EQ:Hormander-estimate-2}) is that it is given in terms of the 
symbol $\sigma_{A}$ of $A$. However, similar to Remark \ref{REM:AR-H}, 
also in the case of a compact Lie group $G$,  Theorem \ref{THM:upper-bound} 
(or rather its version for left-invariant operators in Remark \ref{REM:AR-right})
implies the estimate \eqref{EQ:CG-FM-upper}.
This follows from the following result relating the noncommutative Lorentz norm to the global symbol
of invariant operators (and thus also to representations of $G$)
in the context of compact Lie groups:
\begin{prop} 
\label{PROP:comparison}
Let $1<p\leq 2 \leq q <\infty$ and let $p\not=q$ and $\frac1r=\frac1p-\frac1q$. Suppose $G$ is a compact Lie group and $A$ is a (left) Fourier multiplier on $G$. Then we have
\begin{equation}
\label{EQ:comparison}
\|A\|_{L^{r,\infty}(VN_R(G))}
\leq
\sup_{s>0}
s
\left(
\sum\limits_{\substack{\xi\in\Gh \\ \|\sigma_A(\xi)\|_{\mathcal L(\mathcal H_{\xi})}\geq s}}d^2_{\xi}
\right)^{\frac1p-\frac1q},
\end{equation}
where $\sigma_A(\xi)=\xi^*(g)A\xi(g)\big|_{g=e}\in 
\mathcal L(\mathcal H_{\xi}) \simeq {\mathbb C}^{d_{\xi}\times d_{\xi}}$ 
is the global symbol of $A$.
\end{prop}


\begin{thebibliography}{ANR15}
%
%
%
\bibitem[ANR15]{ANR2015}
Akylzhanov, R., Nursultanov, E., Ruzhansky, M., 2015. Hardy-{L}ittlewood,
  {H}ausdorff-{Y}oung-{P}aley inequalities, and ${L}^p$-${L}^q$ multipliers on
  compact homogeneous manifolds. arXiv:1504.07043.
  
  \bibitem[ANR16]{ANR2016}
Akylzhanov, R., Nursultanov, E., Ruzhansky, M., 2016. 
Hardy-Littlewood inequalities and Fourier multipliers on SU(2). 
to appear in Studia Math.

\bibitem[AR16]{AR2016}
Akylzhanov, R., Ruzhansky, M., 2016. 
  Hausdorff-Young-Paley inequalities and $L^p$-$L^q$ Fourier multipliers on locally compact groups. 
  arXiv:1510.06321.




\bibitem[CW71]{coifman+weiss_lnm}
Coifman, R.~R., Weiss, G., 1971. Analyse harmonique non-commutative sur
  certains espaces homog{\`e}nes.
  Springer-Verlag, Berlin.

%

\bibitem[Di81]{VNA-Dixmier-1981}
Dixmier, J., 1981. Von Neumann algebras. Amsterdam ; New York : North-Holland
  Pub. Co.





\bibitem[FK86]{FackKosaki1986}
Fack T., Kosaki, H., 1986. Generalized $s$-numbers of $\tau$-measurable
  operators. Pacific J. Math. 123~(2), 269--300.

\bibitem[FR14]{Fischer-Ruzhansky:FM-graded}
Fischer, V., Ruzhansky, M., 2014. Fourier multipliers on graded {L}ie groups.
  arXiv:1411.6950.
  
\bibitem[FR16]{Fischer-Ruzhansky:book}
Fischer, V., Ruzhansky, M., 2016. Quantization on nilpotent {L}ie groups.
  Progress in Math., Vol. 314, Birkh{\"a}user.



\bibitem[Ho60]{Hormander:invariant-LP-Acta-1960}
H{\"{o}}rmander, L., 1960. Estimates for translation invariant operators in
  {$L^{p}$}\ spaces. Acta Math. 104, 93--140.
  
\bibitem[HK14]{HK2015}
A.~Hassannezhad and G.~Kokarev.
\newblock Sub-Laplacian eigenvalue bounds on sub-Riemannian manifolds.
\newblock Ann. Scuola Norm. Sup. Pisa Cl. Sci., to appear.


\bibitem[Ko81]{Kosaki1981}
Kosaki, H., 1981. Non-commutative {L}orentz spaces associated with a semifinite
  {V}on {N}eumann algebra and applications. Proc. Japan Acad. Ser. A Math. Sci.
  57~(6), 303--306.

\bibitem[Ku58]{Kunze1958}
Kunze, R.~A., 1958. {$L_{p}$} {F}ourier transforms on locally compact
  unimodular groups. Trans. Amer. Math. Soc. 89, 519--540.

\bibitem[MR15]{MR2015}
Mantoiu, M., Ruzhansky, M., 2015. Pseudo-differential operators, {W}igner
  transform and {W}eyl systems on type {I} locally compact groups.
  http://arxiv.org/pdf/1506.05854v1.pdf.
%
%

\bibitem[Met76]{Metivier1976}
G.~M{\^e}tivier.
\newblock Fonction spectrale et valeurs propres d'une classe d'op\'erateurs non
  elliptiques.
\newblock {\em Comm. Partial Differential Equations}, 1(5):467--519, 1976.
\bibitem[RT10]{RT}
Ruzhansky, M., Turunen, V., 2010. Pseudo-differential operators and symmetries.
  Birkh{\"a}user Verlag, Basel.

\bibitem[RT13]{Ruzhansky+Turunen-IMRN}
Ruzhansky, M., Turunen, V., 2013. Global quantization of pseudo-differential
  operators on compact {L}ie groups, {$\rm SU(2)$}, 3-sphere, and homogeneous
  spaces. Int. Math. Res. Not. IMRN~(11), 2439--2496.

\bibitem[RW15]{RuWi2015}
Ruzhansky, M., Wirth, J., 2015. {$L^p$} {F}ourier multipliers on compact {L}ie
  groups. Math. Z. 280~(3-4), 621--642.

\bibitem[Se53]{Segal1953}
Segal, I., 1953. A non-commutative extension of abstract integration. Ann. 
  Math. 57~(3), 401--457.

\bibitem[Se50]{Segal1950}
Segal, I.~E., 1950. An extension of {P}lancherel's formula to separable
  unimodular groups. Ann. Math. 52, 272--292.




\end{thebibliography}
\end{document}